\title{Making a $K_4$-free graph bipartite}
\author{ Benny Sudakov
\thanks{
Department of Mathematics, Princeton University,
Princeton, NJ 08544. E-mail: bsudakov@math.princeton.edu.
Research supported in part by NSF CAREER award DMS-0546523, NSF grant
DMS-0355497, USA-Israeli BSF grant, and by an Alfred P. Sloan fellowship.}	
}
\date{}
\documentclass [11pt]{article}
\usepackage{latexsym}
\usepackage{amsfonts}
\usepackage{amsmath}
\newtheorem{theo}{Theorem}[section]

\newtheorem{lemma}[theo]{Lemma}
\newtheorem{coro}[theo]{Corollary}

\def\blackslug{\hbox{\kern1pt\vrule height6pt width4pt
depth1pt\kern1pt}}
\def\qed{\penalty 500\hbox{\quad\blackslug}\ifmmode\else\par
\vskip4.5pt plus3pt minus2pt\fi}

\def\t{\triangle}

\begin{document}
\maketitle
\begin{abstract}
We show that every $K_4$-free graph $G$ with $n$ vertices
can be made bipartite by deleting at most $n^2/9$ edges. Moreover,
the only extremal graph which requires deletion of that many edges is 
a complete 3-partite graph with parts of size $n/3$.
This proves an old conjecture of P. Erd\H{o}s.
\end{abstract}

\section{Introduction}
The well-known Max Cut problem asks for the largest bipartite subgraph of a graph $G$. This 
problem has been the subject of extensive research, both from the 
algorithmic perspective
in computer science and the extremal perspective in combinatorics.
Let $n$ be the number of vertices and $e$ be the number edges of $G$ and let
$b(G)$ denote the size of the largest bipartite subgraph of $G$. The extremal part of Max 
Cut problem  asks to estimate $b(G)$ as a function of $n$ and $e$. This question was first 
raised  almost forty years ago by P. Erd\H{o}s \cite{E1} and attracted a lot of attention 
since then (see, e.g., \cite{AKS, ABKS, BS, A, Sh, EGS, EFPS, BL, Ed}).

It is well known that every graph $G$ with $e$ edges can be made bipartite 
by deleting at most $e/2$ edges, i.e., $b(G) \geq e/2$. To see this just consider a random 
partition of vertices of $G$ into two parts $V_1, V_2$ and estimate the expected number of 
edges in the cut $(V_1, V_2)$. A complete graph $K_n$ on $n$ vertices shows that the 
constant $1/2$ in the above bound is asymptotically tight. Moreover, this constant can not 
be improved even if we consider restricted families of graphs, e.g., 
graphs that contain no copy of a fixed {\em forbidden} subgraph $H$. 
We call such graphs $H$-free. Indeed, using sparse random graphs one can 
easily construct a graph $G$ with $e$ edges such that it has no short cycles 
but can not be made bipartite by deleting less than $e/2-o(e)$ edges.
Such $G$ is clearly $H$-free for every forbidden graph $H$ which is not a forest.
It is a natural question to estimate the error term  
$b(G)-e/2$ as $G$ ranges over all $H$-free graph with $e$ edges. 
We refer interested reader to \cite{AKS, ABKS, A, Sh}, where such results were obtained for 
various forbidden subgraphs $H$.

In this paper we restrict our attention to {\em dense} ($e=\Omega(n^2)$)
$H$-free graphs for which it is possible to prove stronger bounds for Max Cut.
According to a long-standing conjecture of Erd\H{o}s \cite{E2}, 
every triangle-free graph on $n$ 
vertices can be made bipartite by deleting at most $n^2/25$ edges. This bound, if true, is 
best possible (consider an appropriate blow-up of a $5$-cycle). 
Erd\H{o}s, Faudree, Pach and Spencer proved that for triangle-free $G$ of order $n$ it is 
enough to delete $(1/18-\epsilon)n^2$ edges to make it bipartite. They also verify 
the conjecture for all graphs with at least $n^2/5$ edges. Some extensions of their results 
were further obtained in \cite{EGS}. Nevertheless this intriguing problem remains open.
Erd\H{o}s also asked similar question for $K_4$-free graphs. 
His old conjecture (see e.g., \cite{EFPS}) asserts that it is enough to
delete at most $(1+o(1))n^2/9$ edges to make bipartite any $K_4$-free graph on $n$ vertices.
Here we confirm this in the following strong form.

\begin{theo}
\label{main}
Every $K_4$-free graph $G$ with $n$ vertices 
can be made bipartite by deleting at most $n^2/9$ edges.
Moreover, the only extremal graph which requires deletion of that many edges is
a complete 3-partite graph with parts of size $n/3$.
\end{theo}

\noindent
This result can be used to prove the following asymptotic generalization.

\begin{coro}
\label{easy}
Let $H$ be a fixed graph with chromatic number $\chi(H)=4$.
If $G$ is a graph on $n$ vertices not containing $H$ as a subgraph, then 
we can delete at most $(1+o(1))n^2/9$ edges from $G$ to make it bipartite.
\end{coro}

Another old problem of Erd\H{o}s, that is similar in spirit, is to 
determine the best local density in $K_r$-free graphs for $r\geq 3$ (for more 
information see, e.g., 
\cite{KS1, KS2} and their references). One of Erd\H{o}s' favorite conjectures was 
that any triangle-free graph $G$ on $n$ vertices should contain a set of $n/2$ 
vertices that spans at most $n^2/50$ edges. The blow-up of a $5$-cycle in which 
we replace 
each vertex by an independent set of size $n/5$ and each edge by a complete 
bipartite graph shows that this estimate can not be improved. On the other 
hand, for $r>3$ Chung and Graham \cite{CG} conjectured that Tur\'an graph has the 
best local density for subsets of size $n/2$. In particular, their conjecture 
implies that every $K_4$-free graph on $n$ vertices should contain a set of $n/2$ 
vertices that spans at most $n^2/18$ edges. 

Krivelevich \cite{K} noticed that for regular graphs a
bound in the local density problem implies a bound for
the problem of making the graph bipartite. Indeed,
suppose $n$ is even, $G$ is a $d$-regular $K_4$-free graph on
$n$ vertices and $S$ is a set of $n/2$ vertices.
Then $dn/2 = \sum_{s \in S} d(s) = 2e(S) + e(S,\bar{S})$
and $dn/2 =  \sum_{s \notin S} d(s) = 2e(\bar{S}) +
e(S,\bar{S})$, i.e. $e(S) = e(\bar{S})$. Deleting
the $2e(S)$ edges within $S$ or $\bar{S}$ makes the
graph bipartite, so if we could find $S$ spanning
at most $n^2/18$ edges we would delete at most $n^2/9$
in making $G$ bipartite. Unfortunately, the converse reasoning does
not work. Nevertheless, we believe that the result of 
Theorem \ref{main} provides some supporting evidence 
for conjecture of Chung and Graham.

The rest of this short paper is organized as follows.
The proof of our main theorem appears in the beginning of next section.
Next we show how to obtain Corollary \ref{easy} 
using this theorem together with well known Szemer\'edi's Regularity 
Lemma \cite{S} (see also \cite{KS}). The last section of the paper 
contains 
some concluding remarks and 
open questions.

\noindent
{\bf Notation.}\, We usually write $G=(V,E)$ for a graph $G$ with
vertex set $V=V(G)$ and edge set $E=E(G)$, setting
$n=|V|$ and $e=e(G)=|E(G)|$. If $X \subset V$ is a
subset of the vertex set then $G[X]$ denotes the
restriction of $G$ to $X$, i.e. the graph on
$X$ whose edges are those edges of $G$
with both endpoints in $X$. We will
write $e(X)=e(G[X])$ and 
similarly, we write $e(X,Y)$ for the
number of edges with one endpoint in $X$ and the
other in $Y$. $N(v)$ is
the set of vertices adjacent to a vertex $v$
and $d(v)=|N(v)|$ is the degree of $v$.  For any two vertices $u,v$ we denote by 
$N(u,v)$ the set of common neighbors of $u$ and $v$, i.e., all 
the vertices adjacent to both of them. We will also write $d(u,v)=|N(u,v)|$. 
Finally if three vertices $u,v,$ and $w$ are all adjacent then they form
a triangle in $G$ and we denote this by $\t=\{u,v,w\}$.

\section{Proofs}
\subsection{Main result}
In this subsection we present the proof of our main theorem. We start with 
the following well known fact (see, e.g., \cite{A}), whose short proof  we include here 
for the sake of completeness.

\begin{lemma}
\label{4-partite}
Let $G$ be a $4$-partite graph with $e$ edges. Then 
$G$ contains a bipartite subgraph with at least $2e/3$ edges.  
\end{lemma}

\noindent
{\bf Proof.}\, Let $V_1, \ldots, V_4$ be a partition of vertices of $G$ into four 
independent sets. Partition these sets randomly into two classes, where each class 
contains exactly two of the sets $V_i$. Consider a bipartite subgraph $H$ of $G$ with these 
color classes. For each fixed edge $(u,v)$ of $G$ the 
probability that $u$ and $v$ will lie in the different classes is precisely $(2\cdot 
2)/{4 \choose 2}=2/3$.
Therefore, by linearity of expectation, the expected number of edges in $H$
is $2e/3$, completing the proof.
\hfill $\Box$

\vspace{0.15cm}
Next we need two simple lemmas.

\begin{lemma}
\label{codegree}
Let $G$ be a graph with $e$ edges and $m$ triangles. Then it contains a 
triangle $\{u,v,w\}$ such that 
$$d(u,v)+d(u,w)+d(v,w) \geq \frac{9m}{e}.$$ 
\end{lemma}

\noindent
{\bf Proof.}\, A simple averaging argument, using that $\sum_{(x,y)\in E(G)} d(x,y)=3m$ 
and Cauchy-Schwartz inequality,  shows that
there is a triangle $\{u,v,w\}$ in $G$ with
\begin{eqnarray*}
\hspace{0.4cm}
d(u,v)+d(u,w)+d(v,w) &\geq& \frac{1}{m}
\sum_{\{x,y,z\}=\t}\Big(d(x,y)+d(x,z)+d(y,z)\Big)=
\frac{1}{m}\sum_{(x,y)\in E(G)} d^2(x,y) \\
&\geq&
\frac{e}{m} \left(\frac{\sum_{(x,y)\in E(G)} d(x,y)}{e}\right)^2
=\frac{(3m)^2}{m e}=\frac{9m}{e}. \hspace{3.4cm} \Box
\end{eqnarray*}

\begin{lemma}
\label{maxcut-1}
Let $G$ be a graph on $n$ vertices with $e$ edges and $m$ triangles. Then 
$G$ contains a bipartite subgraph of size at least $4e^2/n^2-6m/n$.
\end{lemma}

\noindent
{\bf Proof.}\, Let $v$ be a vertex of $G$ and let $e_v$ denotes the number of edges 
spanned by the neighborhood $N(v)$. Consider the bipartite subgraph of $G$ whose parts  
are $N(v)$ and its complement $V(G)\setminus N(v)$. It is easy to see that number of edges in 
this subgraph is $\sum_{u \in N(v)} d(u)-2e_v$. 
Thus averaging over all vertices $v$ we have that
\begin{eqnarray}
\label{A}
b(G) &\geq& \frac{1}{n}\sum_v\bigg(\sum_{u \in N(v)}d(u)-2e_v\bigg)
= \frac{1}{n}\sum_v d^2(v)-\frac{2}{n}\sum_v e_v\\
&\geq& \left(\frac{\sum_v d(v)}{n}\right)^2-6m/n=4e^2/n^2-6m/n \nonumber.
\end{eqnarray}
Here we used Cauchy-Schwartz inequality together with identities $\sum_v e_v=3m$, 
$\sum_v 
d(v)=2e$.
\hfill $\Box$

\vspace{0.15cm}
Now we can obtain our first estimate on the Max Cut in $K_4$-free graphs. This result can be 
used  to prove the conjecture for graphs  with $\leq n^2/4$ edges.

\begin{lemma}
\label{maxcut-2}
Let $G$ be a $K_4$-free graph on $n$ vertices with $e$ edges. Then it 
contains a bipartite subgraph of size at least $2e/7+8e^2/(7n^2)$.
\end{lemma}

\noindent
{\bf Proof.}\, Let $v$ be a vertex of $G$ and denote by $e_v$ the number of edges
spanned by the neighborhood of $v$. Consider a subgraph
of $G$ induced by the set $N(v)$. This subgraph $G[N(v)]$
has $d(v)$ vertices, $e_v$ edges and contains no triangles,
since $G$ is $K_4$-free. Therefore by previous lemma (with $m=0$) it has a bipartite 
subgraph $H$ of size at least $4e_v^2/d^2(v)$. Let $(A,B), A \cup B=N(v)$ be the bipartition of 
$H$. Consider a bipartite subgraph $H'$
of $G$ with parts $(A',B')$, where $A \subset A'$,
$B \subset B'$ and we place each
vertex $v \in V(G) \setminus N(v)$ in $A'$ or $B'$
randomly and independently with probability $1/2$.
All edges of $H$ are edges of $H'$, and
each edge incident to a vertex in $V(G) \setminus N(v)$
appears in $H'$ with probability $1/2$. As the number of edges 
incident to vertices $V(G)\setminus N(v)$ is $e-e_v$,
by linearity of expectation, we have
$b(G) \geq \mathbb{E}\big[e(H')\big] \geq (e-e_v)/2+4e_v^2/d^2(v)$.
By averaging over all vertices $v$ 
\begin{equation}
\label{B}
b(G) \geq \frac{1}{2}e+\frac{1}{n}\sum_v\Big(4e_v^2/d^2(v)-e_v/2\Big).
\end{equation}

To finish the proof we take a convex combination of inequalities (\ref{A}) and (\ref{B})
with coefficients $3/7$ and $4/7$ respectively. This gives

\begin{eqnarray*}
b(G) &\geq& \frac{3}{7}\left(\frac{1}{n}\sum_v d^2(v)-\frac{2}{n}\sum_v e_v\right)  +
\frac{4}{7} \left(\frac{1}{2}e+\frac{1}{n}\sum_v\Big(4e_v^2/d^2(v)-e_v/2\Big)\right)\\
&=&\frac{2}{7}e+\frac{1}{7n}\sum_v\Big(3d^2(v)-8e_v +16e_v^2/d^2(v)\Big)\\
&=&\frac{2}{7}e+
\frac{1}{7n}\sum_vd^2(v)\Big(3-8\big(e_v/d^2(v)\big)+16\big(e_v/d^2(v)\big)^2\Big)\\
&\geq& \frac{2}{7}e +\frac{2}{7n}\sum_vd^2(v) 
\geq \frac{2}{7}e+\frac{2}{7} \left(\frac{\sum_vd(v)}{n}\right)^2=\frac{2}{7}e 
+\frac{8}{7}e^2/n^2,
\end{eqnarray*}
where we used that $3-8t+16t^2=(4t-1)^2+2\geq 2$ for all $t$,
$\sum_vd(v)=2e$ and Cauchy-Schwartz inequality.
\hfill $\Box$

\vspace{0.2cm}
\noindent
{\bf Remark.}\, The above result is enough for our purposes, but
one can get a slightly better inequality by taking a convex combination of 
(\ref{A}) and (\ref{B}) with coefficients $1/(1+a)$ and $a/(1+a)$ with 
$a=1.38$.

\begin{lemma}
\label{technical}
Let $f(t)=t/18+\frac{2}{9}\big(5/2-t-1/t\big)^2$. Then
$f(t)\leq 1/9$ for all $t \in [3/2, 2]$  
and equality holds only when $t=2$.
\end{lemma}

\noindent
{\bf Proof.}\, Note that $f(2)=1/9$ and 
$$f(t)-1/9=\frac{4t^4-19t^3+31t^2-20t+4}{18t^2}=
\frac{(t-2)(4t^3-11t^2+9t-2)}{18t^2}.$$
Consider $g(t)=4t^3-11t^2+9t-2$ in the interval $[3/2, 2]$.
The derivative of this function $g'(t)=12t^2-22t+9$ is zero when
$t=\frac{22\pm\sqrt{52}}{24}$, so the largest root of $g'(t)$ is less than $3/2$.
Therefore $g(t)$ is strictly increasing function for $t \geq 3/2$ and 
so  $g(t)>g(3/2)=1/4>0$ for all $t\in [3/2, 2]$. Since $18t^2>0$ and $t-2$ is negative for 
$t<2$ we conclude that $f(t)-1/9<0$ for all $t \in [3/2, 2)$.
\hfill $\Box$

\vspace{0.15cm}
Having finished all the necessary preparations we are now in a position to complete the 
proof of our main result.

\vspace{0.2cm}
\noindent
{\bf Proof of Theorem \ref{main}.}\, It is easy to see that 
complete 3-partite graph with parts of size $n/3$ has
$(n/3)^3=n^3/27$ triangles and that every edge of this graph is contained in 
exactly $n/3$ of them. To make this graph bipartite we need to
destroy all these triangles. Since deletion of one edge can 
destroy at most $n/3$ of them, altogether we need to delete at least
$\frac{n^3/27}{n/3}=n^2/9$ edges. To finish the proof it remains to show that 
deletion of $\leq n^2/9$ edges is sufficient to make every $K_4$-free graph
bipartite.

Let $G$ be a $K_4$-free graph on $n$ vertices with $e$ 
edges. Tur\'an's theorem \cite{T} says that $e \leq n^2/3$, with equality 
only when $G$ is a 
complete 3-partite graph
with parts of size $n/3$. By Lemma \ref{maxcut-2}, we need to delete at most
$e-b(G) \leq 5e/7-8e^2/(7n^2)=\big(\frac{5}{7}(e/n^2)-\frac{8}{7}(e/n^2)^2\big)n^2$ 
edges to make $G$ bipartite. The function $g(t)=5t/7-8t^2/7$ is increasing in the interval
$t \leq 1/4$ and so $g(t) \leq g(1/4)=3/28$. Therefore if $e\leq n^2/4$ we 
can delete at most $3n^2/28<n^2/9$ edges to make $G$ bipartite. 

Next, consider the case when $n^2/4\leq e \leq n^2/3$ and let $m$ be the number of triangles in 
$G$. By Lemma \ref{maxcut-1}, we can delete at most $e-b(G) \leq e-\big(4e^2/n^2-6m/n\big)$ 
edges to make $G$ bipartite. So we can assume that 
$e-4e^2/n^2+6m/n \geq n^2/9$ or we are done. Then the number of triangles in $G$ satisfies
$m\geq \frac{n}{6}\big(n^2/9+4e^2/n^2-e\big)$ and Lemma \ref{codegree} implies that $G$ 
contains a triangle $\t=\{u,v,w\}$ with 
$$d(u,v)+d(u,w)+d(v,w) \geq \frac{9m}{e}\geq 6e/n + n^3/(6e) -3n/2.$$
Let $V_1=N(u,v), V_2=N(u,w)$, $V_3=N(v,w)$ and let
$X=V(G)\setminus(\cup_{i=1}^3 V_i)$. Since $G$ is $K_4$-free and $(u,v), (u,w),
(v,w)$ are edges of $G$ we have that sets $V_i, 1\leq i\leq 3$ are independent
and disjoint. Consider a 4-partite subgraph $G'$ of $G$ with parts $V_1, V_2, V_3$ and $X$.
This graph has $e(G')=e-e(X)$ edges where $e(X)$ is the number of edges spanned
by $X$. By Tur\'an's theorem $e(X) \leq |X|^2/3$ and we also know that
$$|X|=n-\sum_i |V_i|=n-\big(d(u,v)+d(u,w)+d(v,w)\big) \leq
5n/2-6e/n-n^3/(6e).$$

Since $G'$ is 4-partite we can now use Lemma \ref{4-partite} to 
deduce that $b(G)\geq b(G') \geq 2e(G')/3 =\frac{2}{3}\big(e-e(X)\big).$ 
Therefore the number of edges we need to delete to make
$G$ bipartite is bounded by
\begin{eqnarray*}
e-b(G) &\leq& e-2\big(e-e(X)\big)/3=e/3+2e(X)/3 \leq e/3+2|X|^2/9\\
&\leq& e/3+\frac{2}{9} \Big(5n/2-6e/n-n^3/(6e)\Big)^2\\
&=&\left(\frac{1}{18}(6e/n^2)+\frac{2}{9}
\Big(5/2-6e/n^2-(6e/n^2)^{-1}\Big)^2\right)n^2\\
&=&f\big(6e/n^2\big) \cdot n^2,
\end{eqnarray*}
where $f(t)=t/18+\frac{2}{9}\big(5/2-t-1/t\big)^2$.
As $n^2/4 \leq e \leq n^2/3$ we have that $3/2 \leq t=6e/n^2\leq 2$. Then,
by Lemma \ref{technical}, $f(6e/n^2) \leq 1/9$ with equality only if 
$e=n^2/3$. This shows that we can delete at most $n^2/9$ edges to make $G$
bipartite and we need to delete that many edges only when $e(G)=n^2/3$, 
i.e., $G$ is a complete 3-partite graph
with parts of size $n/3$. .  
\hfill $\Box$

\subsection{Forbidding fixed $4$-chromatic subgraph}
In this short subsection we show how to use Theorem \ref{main} to deduce 
a similar statement about graphs  with any fixed 
forbidden $4$-chromatic subgraph. 
The proof is a standard application of Szemer\'edi's Regularity Lemma and 
we refer the interested reader to the excellent survey of Koml\'os and 
Simonovits \cite{KS}, which discusses various results proved by this powerful 
tool.

We start with a few definitions, most of which follow \cite{KS}.
Let $G=(V,E)$ be a graph, and let $A$ and $B$ be two disjoint 
subsets of $V(G)$. If $A$ and $B$ are non-empty,
define the {\em density of edges} between $A$ and $B$ by 
$ d(A,B) = \frac{e(A,B)}{|A||B|}$.
For $\epsilon>0$ the pair $(A,B)$ is called 
{\em $\epsilon$-regular} if 
for every $X \subset A$ and $Y \subset B$ satisfying
$|X|>\epsilon |A|$ and $ |Y|>\epsilon |B|$
we have $|d(X,Y)-d(A,B)| < \epsilon$.
An {\em equitable partition} of a set $V$ is a partition of $V$ into 
pairwise disjoint classes $V_1,\cdots,V_k$ of almost equal size, 
i.e., $\big| |V_i|-|V_j| \big| \leq 1$ for all $i,j$.
An equitable partition of the set of vertices $V$ of $G$ into 
the classes $V_1,\cdots,V_k$ is called 
{\em $\epsilon$-regular} if $|V_i| \leq \epsilon |V|$ for every $i$
and all but at most $\epsilon k^2$ of the pairs $(V_i,V_j)$ are
$\epsilon$-regular. The above partition is called 
{\em totally $\epsilon$-regular} if all the pairs $(V_i,V_j)$ are 
$\epsilon$-regular.
The following celebrated lemma was proved by 
Szemer\'edi in \cite{S}.

\begin{lemma}
\label{szemeredi}
For every $\epsilon>0$  there is an integer $M(\epsilon)$
such that every graph of order $n>M(\epsilon)$
has an $\epsilon$-regular partition into $k$ classes,
where $k \leq M(\epsilon)$.
\end{lemma}

In order to apply the Regularity Lemma we need to show the existence of a 
complete multipartite subgraph in graphs with a totally $\epsilon$-regular 
partition. This is established in the following lemma which is a 
special case of a well-known result, see, e.g., \cite{KS}.

\begin{lemma}
\label{key}
For every $\delta>0$ and integer $t$ there exist an $0<\epsilon=
\epsilon(\delta,t)$ and $n_0=n_0(\delta,t)$ with the following 
property. 
If $G$ is a graph of order $n>n_0$ and $(V_1,\cdots,V_4)$
is a totally $\epsilon$-regular partition of vertices of $G$ such that 
$d(V_{i},V_{j})\geq \delta$ for all $i<j$, then $G$ contains 
a complete $4$-partite subgraph $K_4(t)$ with parts of size $t$.
\end{lemma}

\noindent 
{\bf Proof of Corollary \ref{easy}.}\,
Let $H$ be a fixed $4$-chromatic graph of order $t$ and
let $G$ be a graph on $n$ vertices not containing $H$ as a subgraph.
Suppose $\delta>0$ and let 
$\epsilon=\min\big(\delta, \epsilon(\delta,t)\big)$, where 
$\epsilon(\delta,t)$ is defined in the previous statement.
Then, by Lemma \ref{szemeredi}, for sufficiently large $n$ there exists an 
$\epsilon$-regular partition $(V_1,\cdots,V_k)$ of vertices of $G$.

Consider a new graph $G'$ on the vertices $\{1, \ldots, k\}$ in which 
$(i,j)$ is an edge iff $(V_i,V_j)$ is an $\epsilon$-regular pair with 
density at least $\delta$. We claim that $G'$ contains no $K_4$. 
Indeed, any such clique in $G'$  corresponds to $4$ parts in the partition of $G$ such that any 
pair of them is $\epsilon$-regular and has density at least $\delta$. This contradicts our 
assumption on $G$, since by Lemma \ref{key}, the union of  these parts will contain a copy of 
complete $4$-partite graph $K_4(t)$ which clearly contains $H$. 

By applying Theorem \ref{main} to graph $G'$, we deduce that there is a set
$D$ of at most $k^2/9$ edges of $G'$ whose deletion makes it bipartite.
Now delete all the edges of $G$ between the pairs $(V_i,V_j)$ with $(i,j) \in D$.
Delete also the  edges of $G$ that lie
within classes of the partition, or that belong to a non-regular pair, or that join a pair of 
classes of density less than $\delta$. It is easy to see that the remaining graph is bipartite
and the number of edges we deleted is at most
$$\hspace{3.2cm}
(k^2/9)(n/k)^2+\epsilon n^2+\delta n^2 \leq \big(1/9+2\delta\big) n^2=(1+o(1))n^2/9. 
\hspace{3.2cm} 
\Box$$

\section{Concluding remarks}
How many edges do we need to delete to make a $K_r$-free graph $G$ of 
order $n$
bipartite? For $r=3, 4$ this was asked long time ago by P. Erd\H{o}s. For triangle-free 
graphs he conjectured that deletion of $n^2/25$ edges is always enough and that 
extremal example is a blow-up of a $5$-cycle. In this paper we answered the question for
$r=4$ and proved that the unique extremal construction in this case is a 
complete $3$-partite graph with equal parts. Our result suggests that a 
complete $(r-1)$-partite graph of order $n$ with equal parts is worst example also 
for all remaining values of $r$. Therefore we believe that 
it is enough to delete at most $\frac{(r-2)^2}{4(r-1)^2}n^2$ edges for 
even $r\geq 5$ and at most $\frac{r-3}{4(r-1)}n^2$ edges for odd $r\geq 5$ 
to make bipartite any $K_r$-free graph $G$ of order $n$. 
It seems that some of the ideas presented here can be useful to make a 
progress on this problem for even $r$.

\vspace{0.4cm}
\noindent
{\bf Acknowledgment.} I would like to thank J\'ozsef Balogh
and Peter Keevash for interesting discussions on the early stages of this project.

\end{document}